\documentclass[12pt]{article}
\usepackage{graphicx}
\usepackage{amssymb}
\usepackage{amsmath}
\usepackage{amsfonts}
\usepackage{amsthm}
\usepackage[english]{babel}
\usepackage[latin1,ansinew]{inputenc}
\usepackage{a4wide}
\usepackage{color}
\newcommand{\be}{\begin{eqnarray}}
\newcommand{\ee}{\end{eqnarray}}
\newtheorem{theo}{Theorem}

\newcommand{\R}{\mathbb R}

\newcommand{\I}{\mathbb I}

\newcommand{\eps}{\epsilon}

\begin{document}
\date{December 7, 2021}

\title{Remarks on Ruggeri's First Model of Dissipative 
\\ Fluid Dynamics as a Symmetric Hyperbolic System} 

\author{Heinrich Freist\"uhler}
\maketitle
\vspace{-.5cm}
The first appealing formulation for the dynamics of viscous heat-conductive fluids as
a symmetric hyperbolic system of balance laws was given by Ruggeri in his pioneering 1983
paper \cite{R83}. Using   
$$
\Upsilon=(\Upsilon_1,\Upsilon_2,\Upsilon_3,\Upsilon_4,\Upsilon_5,\Upsilon_6)=
(\tilde\psi,\tilde u,\tilde\theta,\tilde\Sigma,\tilde\sigma,\tilde q) 
$$
with 
\be\label{GVnonrelExtended} 
\tilde\psi\equiv\frac{g-u^2/2}{\theta},
\ \tilde u\equiv\frac u\theta,
\ \tilde\theta=-\frac1\theta,
\ \tilde\Sigma\equiv\frac\Sigma\theta,
\ \tilde\sigma\equiv\frac\sigma\theta,  
\ \tilde q \equiv\frac{q}{\theta^2},
\ee
as state variable, this formulation reads
\be\label{RSHS}
\frac\partial{\partial t}\left(\frac{\partial X^0(\Upsilon)}{\partial \Upsilon_i}\right)
+
\sum_{i=1}^3\frac\partial{\partial x^j}\left(\frac{\partial X^j(\Upsilon)}{\partial \Upsilon_i}\right)
=
I^i(\Upsilon),\quad i=1,...,6,
\ee
where $X^0$ and $X^1,X^2,X^3$ are the Godunov-Boillat potentials \cite{Go,Bo}
$$
X^0=\frac p\theta,\quad 
X^j=
\left(p\I+\theta(\tilde\Sigma +\tilde\sigma\I)\right) \tilde u_j+\theta\tilde q_j
$$
with
$$
p=
\hat p
\left(
\theta,\tilde\psi+\frac12\theta\tilde u^2
+\frac12\tau_1\tilde\Sigma:\tilde\Sigma
+\frac12\tau_2\tilde\sigma^2
+\frac12\tau_0\tilde q^2
\right)
$$
while the source $I(\Upsilon)$ is given by 
$
I^1(\Upsilon)=0_1,
I^2(\Upsilon)=0_3,
I^3(\Upsilon)=0_1
$ 
and
$$
I^4(\Upsilon)=-\frac\Sigma{2\eta},\
I^5(\Upsilon)=-\frac\sigma{3\zeta},\
I^6(\Upsilon)=-\frac q\chi.
$$
Here, $u\in\R^3,\theta,g,-(\Sigma+\sigma\I)$ (with $\Sigma\in\R^{3\times 3}$ tracefree and symmetric),
and $q\in \R^3$ denote velocity, temperature, chemical potential, viscous stress, and heat flux, while 
the fluid is specified 
through an equation of state
$$
p=\hat p(\theta,\psi)
$$ 
with the extended thermodynamical potential 
$$
\psi= 
\tilde\psi
+\frac12\theta\tilde u^2
+\frac12\tau_1\tilde\Sigma:\tilde\Sigma
+\frac12\tau_2\tilde\sigma^2
+\frac12\tau_0\tilde q^2,
$$    
and constants $\tau_0,\tau_1,\tau_2>0$ .
While Ruggeri and M\"uller later, in their extensive work leading to the fundamental theory of 
\emph{Rational Extended Thermodynamics} (RET), went on to more refined formulations based on 
\emph{main fields} (cf.\ \cite{RS}) that are different from \eqref{GVnonrelExtended} (cf.\ \cite{MR}), it seems  
that the abovedescribed 
`Ruggeri's first model' still deserves attention at least from mathematical points of view, as it 
may serve as a prototype regarding the latter. 

The purpose of this note is to point out two aspects. The first one is related to the fact that 
if one replaces $\tau_0,\tau_1$ and $\tau_2$ by $0$, equations \eqref{RSHS} become the Navier-Stokes-Fourier
(NSF) equations (cf.\ \cite{R83}); we fix $\bar\tau_0,\bar\tau_1,\bar\tau_2>0$ and denote for $\epsilon\ge 0$
system \eqref{RSHS} with 
$$
(\tau_0,\tau_1,\tau_2)=\epsilon(\bar\tau_0,\bar\tau_1,\bar\tau_2)
$$
as \eqref{RSHS}$_\epsilon$.
\begin{theo}
If for some $T>0$ the NSF equations \eqref{RSHS}$_0$ with given initial data have a solution with time 
sections, for $0\le t\le T$, in an appropriate Sobolev space, then so does Ruggeri's model 
\eqref{RSHS}$_\eps$, 
if $\epsilon>0$ is sufficiently small. The solutions to the latter tend, for $\epsilon\searrow 0$, to the NSF solution. 
\end{theo}
\begin{theo}
For any fixed $\tau_0,\tau_1,\tau_2>0$, 
and any homogenenous reference state with $\theta,\psi,u$ constant and $\Sigma=\sigma=q=0$,  
data that in an appopriate Sobolev norm deviate sufficiently little from the reference state 
lead to solutions of \eqref{RSHS} that exist globally in time and decay for $t\to\infty$ to the 
reference state.
\end{theo}
Theorem 1 is modelled on an analogous result that was established by Yong in his admirable paper 
\cite{Y14} on a similar model for barotropic 
fluids, of which Ruggeri's system \eqref{RSHS} is a  non-barotropic version.\footnote{Ruggeri's system
is moreover Galilei invariant. Cf.\ \cite{GIVYM} for a Galilei invariant version of Yong's model.} Indeed,
it is shown in complete analogy to Yong's proof.

The assertion of Theorem 2 is analogous to a result Ruggeri obtained in 
\cite{R04} on the RET formulation of dissipative relativistic fluid dynamics \cite{MR}. 
The proof of Theorem 2 is given below. 

{\bf Remark.} Analoga of Theorems 1 and 2 hold for corresponding balance-law models of barotropic fluids,
and for counterparts of both the latter and the original Ruggeri model in the relativistic contexts as 
proposed in \cite{GIVYM}; they are proved in the same way. 
 
\emph{Proof of Theorem 2.} We first compute the derivatives of the potentials  $X^0$ and $X=(X^1,X^2,X^3)$
with respect to the Godunov variables. Using $\hat p_\psi=\rho\theta$ and $\theta\hat p_\theta=\rho e+p$, 
we find
\begin{align*}
\frac{\partial X^0}{\partial\tilde\psi}&=\frac{\hat p_\psi}\theta=\rho  
\\
\frac{\partial X^0}{\partial\tilde u}&= \hat p_\psi\tilde u=\rho u 
\\
\frac{\partial X^0}{\partial\tilde\theta}&
= -p+\theta\frac{\partial p}{\partial\theta}
=-\hat p+\theta\hat p_\theta+\frac12 \hat p_\psi \theta\tilde u^2 
=\rho e+\frac12\rho u^2
\\
\frac{\partial X^0}{\partial\tilde\Sigma}&
=\frac1\theta\hat p_\psi\tau_1\tilde\Sigma
=\tau_1\rho\frac\Sigma\theta   
\\
\frac{\partial X^0}{\partial\tilde\sigma}&
=\frac1\theta\hat p_\psi\tau_2\tilde\sigma
=\tau_1\rho\frac\sigma\theta   
\\
\frac{\partial X^0}{\partial\tilde q}&
=\frac1\theta\hat p_\psi\tau_0\tilde q
=\tau_0\rho\frac q{\theta^2}
\end{align*}
and
\begin{align*}
\frac{\partial X}{\partial\tilde\psi}&=\hat p_\psi\tilde u=\rho u     \\
\frac{\partial X}{\partial\tilde u}&
=
(\hat p+\hat p_\psi\theta \tilde u^2)\I+\theta(\tilde\Sigma+\tilde\sigma\I)
= 
(p+\rho u^2)I+(\Sigma+\sigma\I) \\
\frac{\partial X}{\partial\tilde\theta}&
=
\theta^2\left(\left(\left(\hat p_\theta+\frac12\hat p_\psi \tilde u^2\right)I
+(\tilde\Sigma+\tilde\sigma\I)\right)\tilde u
+\tilde q\right)=
\left(\left({\rho e+p}+\frac12{\rho u^2}\right)I+(\Sigma+\sigma\I)\right) u + q \\
\frac{\partial X}{\partial\tilde\Sigma}&=
(\tau_1\hat p_\psi\tilde\Sigma+\theta\I)\tilde u  
=\left(\tau_1\rho\frac\Sigma\theta+\I\right) u
\\
\frac{\partial X}{\partial\tilde\sigma}&=
(\tau_1\hat p_\psi\tilde\sigma+\theta)\tilde u  
=\left(\tau_1\rho\frac\sigma\theta +1\right) u
\\
\frac{\partial X}{\partial\tilde q}&=    
\tau_0\hat p_\psi\tilde q\tilde u+\theta
=\tau_0\rho\frac q{\theta^2}u+\theta.
\end{align*}
Then, for $u=\tau=q=0$, the Hessians of $X^0$ and, w.\ l.\ o.\ g., $X^1$ are   
$$
D^2X^0=\begin{pmatrix}
       \theta^{-1}\hat p_{\psi\psi}&0&-\hat p_\psi+\theta \hat p_{\theta\psi}&0&0&0\\
       0&\hat p_\psi&0&0&0&0\\
       -\hat p_\psi+\theta \hat p_{\theta\psi}&0&\theta\hat p_{\theta\theta}&0&0&0\\
       0&0&0&\tau_1\theta^{-1} \hat p_\psi&0&0     \\
       0&0&0&0&\tau_2\theta^{-1} \hat p_\psi&0     \\
       0&0&0&0&0&\tau_0\theta^{-1} \hat p_\psi      
       \end{pmatrix}
$$
and 
$$
D^2X^1=\begin{pmatrix}
       0&\hat p_\psi&0&0&0&0\\
       \hat p_\psi&0&\theta^2 p_\theta&\theta&\theta&0\\
       0&\theta^2\hat p_\theta&0&0&0&\theta^2\\ 
       0&\theta&0&0&0&0\\
       0&\theta&0&0&0&0\\
       0&0&\theta^2&0&0&0      
       \end{pmatrix},
$$
while the Jacobian of the source is   
$$
DI=\begin{pmatrix}
   0&0&0&0&0&0\\
   0&0&0&0&0&0\\
   0&0&0&0&0&0\\
   0&0&0&-1/(2\eta)&0&0\\
   0&0&0&0&-1/(3\zeta)&0\\
   0&0&0&0&0&-1/\chi
\end{pmatrix}.
$$
According to \cite{R04,BHN} and because of Galilean invariance,
it suffices to check the Kawashima condition \cite{K}:
\be\label{K}
\text{No eigenvector of $D^2X^1$ with respect to $D^2X^0$ lies in ker$(DI)$.} 
\ee
Now, ker$(DI)$ consists exactly of all vectors of the form $v=(v_1,v_2,v_3,0,0,0)^\top$ and if for any 
such $v$ and any $\lambda\in\R$, 
$$
(-\lambda D^2X^0+D^2X^1)v=0,  
$$
then necessarily $v_2=v_3=0$ and thus  
$$
\text{$v=0$ or $p_\psi=0$.}
$$ 
As $p_\psi=\rho\theta>0$, we have established \eqref{K}. 
$\hspace*{8.7cm}\Box$

\end{document}